\documentclass[11pt]{article}

\usepackage{babel}
\usepackage{amsmath,amsthm}
\usepackage[T1]{fontenc}
\usepackage[utf8]{inputenc}
\usepackage{amsfonts}
\usepackage{amssymb}
\usepackage{graphicx}
\usepackage{amscd}

\usepackage{amscd}
\usepackage[all]{xy}
\usepackage[usenames, dvipsnames]{color}
\usepackage{fancyhdr}
\usepackage{makeidx}
\makeindex

\theoremstyle{definition} 
 \newtheorem{thm}{Theorem}[section]
 \theoremstyle{definition}
 
 \newtheorem{cor}[thm]{Corollary}
 \newtheorem{lem}[thm]{Lemma}
 \newtheorem*{lem*}{Lemma}
 \newtheorem*{defn*}{Definition}
 \newtheorem{prop}[thm]{Proposition}
 \theoremstyle{definition}
 \newtheorem{defn}[thm]{Definition}
 \theoremstyle{remark}
 \newtheorem{rem}[thm]{Remark}
 \newtheorem*{rem*}{Remark}
  \newtheorem*{thm*}{Theorem}
 
\newtheorem{nota}[thm]{Notation}
 \newtheorem{ejem}[thm]{Example}

\numberwithin{equation}{section}
\newcommand{\OO}{{\mathcal O}}

\newcommand{\A}{{\mathcal A}}

\newcommand{\U}{{\mathcal U}}

\newcommand{\PP}{{\mathcal P}}

\newcommand{\X}{{\mathcal X}}
\newcommand{\op}{{\mathrm{op}}}
\newcommand{\qc}{{\mathrm{qc}}}

\newcommand{\PPP}{{\mathbf P}}

\newcommand{\D}{{\mathfrak D}}

\newcommand{\V}{{\mathcal V}}

\newcommand{\F}{{\mathcal F}}
\newcommand{\G}{{\mathcal G}}
\newcommand{\C}{{\mathfrak C}}

\newcommand{\x}{{\mathbf{x}}}

\DeclareMathOperator{\Qcoh}{{\mathbf{Qcoh}}}

\DeclareMathOperator{\Spec}{\mathrm{Spec}}

\DeclareMathOperator{\Rev}{\textbf{RÉt}}

\DeclareMathOperator{\Fib}{\mathrm{Fib}}
\DeclareMathOperator{\tlim}{2\text{-lim}}

\DeclareMathOperator{\tcolim}{2\text{-colim}}

\DeclareMathOperator{\colim}{\mathrm{colim}}
\DeclareMathOperator{\data}{{\text{-}\mathbf{data}}}

\newcommand{\p}{{\mathfrak p}}

\newcommand{\enumera}{\begin{enumerate}}
\newcommand{\eenumera}{\end{enumerate}}
\DeclareMathOperator{\pos}{{\mathbf{pos}}}

\DeclareMathOperator{\Hom}{\mathrm{Hom}}
\newcommand{\Cyl}{{\mathrm{Cyl}}}

\begin{document}
\title{Lax colimits of posets with structure sheaves: applications to descent}
\author{Javier Sánchez González\\ \small{Universidad de Salamanca, Department of Mathematics}\\ javier14sg@usal.es}

\maketitle
\begin{abstract}
We consider categories of posets with $\C$-valued structure sheaves  for any category $\C$ and see how they possess poset-indexed lax colimits that are both easy to describe and "weakly equivalent" to their ordinary colimits in a certain sense. We employ this construction to study descent problems on schematic spaces---a particular scheme-like kind of ringed poset---, proving a general Seifert-Van Kampen Theorem for their étale fundamental group that recovers and generalizes the homonym result for schemes to the topology of flat monomorphisms. The techniques are general enough to consider their applications in many other frameworks.
\end{abstract}
\section{Introduction}

In the geometric world, categorical colimits can be thought of as a general way of expressing gluing of spaces. For example, any scheme is the colimit, in the category of locally ringed spaces, of the components of any of its affine coverings. Furthermore, for any reasonable scheme, one might assume that these colimits are indexed by posets: the \textit{nerve} of the corresponding affine covering with its redundancies removed. This idea is simply a generalization to locally ringed spaces of the construction of \textit{finite models} of topological spaces, an old technique of McCord to study homotopy types \cite{mccord}.

To better understand the situation, one may consider these \textit{recollement} data of schemes as functors from some poset $X$, understood as a category,  which we may assume finite under mild compactness hypothesis on the original scheme $S$; to the category of commutative rings with unit, i.e. $X\to \mathbf{CRing}$. It is also very classical that the category of such functors coincides with the category of sheaves of rings on the poset $X$, understood as a topological space, so we may study these collections of data as (non-locally) ringed spaces themselves. A major advantage of this point of view is that many sheaf-theoretic notions admit a simple description, like quasi-coherent modules or the \v{C}ech resolution of an abelian sheaf; and in many cases they coincide with their scheme-theoretic analogues on $S$. From these observations, Sancho first axiomatized in \cite{fernando schemes} a---non-full---subcategory of finite ringed posets that behaves suitably well with respect to quasi-coherent modules and that, not only contains all finite models of (quasi-compact and quasi-separated) schemes, but also that non-trivially generalizes them: the category of (finite) schematic spaces, $\mathbf{SchFin}$. After spending some time with these objects, it is easy to convince oneself of their geometric interest compared to similarly-purposed constructions such as those of simplicial schemes. A brief summary of this non-standard approach will be provided in Section \ref{section schematic}.

Taking a step back from the previous discussion, we also note that computing colimits in the category of schemes---or simply determining if they exist---is a very general and difficult problem: even when existence is guaranteed---for example, colimits via open immersions---, obtaining explicit expressions for them is no easy task. The situation is not much better if we attempt to compute colimits of schematic spaces as described in the previous paragraph; however, due to their combinatorial nature, there is an alternative: given a \textit{poset-indexed} "datum" $\mathcal{U}\colon P\to \mathbf{SchFin}^\op$, we may construct a space, which we call the \textit{cylinder of $\mathcal{U}$}, by simply turning the set $\coprod_{p\in P}|\mathcal{U}(p)|$ into a poset with the structure inherited from the underlying posets of each $|\mathcal{U}(p)|$ and the transition morphisms $|\U(p)|\to |\U(q)|$ for $q\leq p$. If we endow the resulting poset with the structure sheaf---which is just a functor---induced by the structure sheaves of each $\U(p)$, we obtain a ringed poset---which will be another schematic space under certain conditions---that \textit{represents}---will be \textit{weakly equivalent} to---the desired colimit in a precise sense and which has been computed without performing any complicated categorical operations on either commutative rings or posets. 

In this paper, we will study this \textit{cylinder} space in its most general formulation, replacing $\mathbf{CRing}$ for any category $\C$. The category of \textit{finite posets with a $\C$-values structure sheaf} will be constructed and called \textit{the category of $\C$-data}, denoted $\C\data$. If $\C=\pos$ is the category of posets, it turns out that the cylinder space is just an incarnation of a \textit{poset-indexed lax colimit}, understanding $\pos$ as a strict 2-category in the natural way. We prove this and see how it generalizes to posets admiting structure sheaves with values in $\C$, detailing the different 2-categorical structures that one may consider and how they interact with each other. For this purpose, we will need to consider $\D\data$ with $\D$ being a strict 2-category, which also proves useful in other applications. Our language of choice will keep things analogue to the category of ringed posets that we use as a reference.

These constructions, of $\C\data$ and lax colimits on it, are general enough to model many descent problems, which appear in a natural and functorial way. For the sake of keeping the discussion focused, we will specialize everything to the previously-discussed schematic case: we will characterize when a cylinder of schematic spaces remains schematic and will apply such characterization to give a very general descent Theorem for \textit{data} on schematic spaces, which even admits a topos-theoretic interpretation---only sketched here due to space limitations---. We will see some examples, with the main one being the Seifert-Van Kampen Theorem for the \'etale fundamental group of schematic spaces, as constructed in \cite{paper grupo etale}. It is worth noting that a variation of the homonym result for schemes and their topology of flat monomorphisms---rather than Zariski or étale---follows purely from the formal descent result developed in this paper and the classical case, which exemplifies how these techniques automatically extend any "reasonable" Zariski-local statement to the aforementioned topology.

\section{Motivation: schematic spaces}\label{section schematic}
Let us give a rather informal introduction to the objects of study for the applications, which were the original motivation to develop the theory of $\C\data$ that we will introduce in the following sections. In prose: schematic spaces arise as the largest subcategory of ringed finite posets that behaves "like quasi-compact and quasi-separated (qc-qs)" schemes with respect to categories of quasi-coherent sheaves. The basic example comes from the construction of finite models of schemes, see \cite{fernando homotopy}, which is a generalization of an earlier topological technique, see \cite{mccord}: given a qc-qs scheme $S$ and a finite covering $\{U_i\}$ of $S$, one may define a poset $X$ as the T$_0$-fication of the topology generated by the covering. Explicitly, if for any $s\in S$ we denote $U^s=\cap_{s\in U_i}U_i$, one sets $X=S/\sim$ with $s\sim s'$ whenever $U^s=U^{s'}$ and $[s]\leq [s']$ if and only if $U^{s'}\subseteq U^s$. The result is a morphism of ringed spaces
\begin{align*}
\pi\colon S\to (X, \pi_*\OO_S).
\end{align*}
If the covering is chosen so that the $U^s$ are affine, $\pi$ induces an adjoint equivalence $(\pi^*\dashv \pi_*)\colon \Qcoh(S)\overset{\sim}{\to}\Qcoh(X)$, so $S$ can be studied from $X$.

Schematic spaces were first introduced in \cite{fernando schemes} and studied in \cite{fernando afines}, \cite{fernando dualidad} or \cite{paper grupo etale}. In \cite{notas de pedro} one can find an extensive compilation of characterizations, and the author of this paper has expanded upon this and exhaustively explored their role in algebraic geometry in his PhD thesis, which includes the contents of this paper and more to come. Morally speaking, we recommend to think of a schematic space $X$ as a model or "structured descent data" of the locally affine locally ringed space $\Spec(X)=\colim_{x\in X}\Spec(\OO_{X, x})$, but the schematic condition forces such discrete incarnation $X$ to be "nice enough" to reflect the geometry of $\Spec(X)$ to some extent. One shows that these $\Spec(X)$ have to be locally affine in the topology of flat monomorphisms of affine schemes and, thus, contain all qc-qs schemes; but schematic spaces do not model, for example, algebraic spaces, for which the associated locally ringed space does not preserve enough useful algebraic information. There are also "geometric" arguments for considering ringed posets as the basis for our combinatorial models over other, more classical, alternatives like simplicial schemes. 

While perhaps not the most enlightening approach, it will be convenient for our purposes to consider the following definitions. Assume for simplicity that all stalk rings of our ringed posets are Noetherian.
\begin{defn}\label{definition schematic}
A finite ringed poset $X$ is a (finite) \textit{schematic space} if
\begin{itemize}
\item For any $x\leq y$, the morphism $r_{xy}\colon \OO_{X, x}\to\OO_{X, y}$ is flat.
\item For any $t\leq x, y$, the morphism
\begin{align*}
\OO_{X, x}\otimes_{\OO_{X, t}}\OO_{X, y}\to \prod_{z\geq x, y}\OO_{X, z}
\end{align*}
is faithfully flat.
\end{itemize}
A morphism $f\colon X\to Y$ between schematic spaces is \textit{schematic} if
\begin{itemize}
\item For any $x\in X$ and $y\geq f(x)$, the morphism
\begin{align*}
\OO_{X, x}\otimes_{\OO_{Y, f(x)}}\OO_{Y, y}\to \prod_{z\in U_x\cap f^{-1}(U_y)}\OO_{X, z}
\end{align*}
induces a surjection between the prime spectra.
\end{itemize}
\end{defn}
\begin{rem}\label{remark flat epic}
A simple descent argument for faithfully flat morphisms shows that $\OO_{X, x}\otimes_{\OO_{X, y}}\OO_{X, y}\simeq \OO_X(U_x\cap U_y)$ for all $t\leq x, y$. If $x=y$, this condition implies that the restriction morphisms $r_{xy}$ of any schematic space are \textit{flat epimorphisms of rings}, hence local isomorphisms.
\end{rem}

Let $\mathbf{SchFin}$ denote the category of schematic spaces and morphisms. All ringed posets and morphisms will be considered schematic unless stated otherwise. The $\Spec$ construction outlined in the introduction of this section defines a functor to the category of locally ringed spaces which is neither full or faithful:
\begin{align*}
\Spec\colon \mathbf{SchFin}\to \mathbf{LRS} & & X\mapsto \Spec(X):=\colim_{x\in X}\Spec(\OO_{X, x}).
\end{align*}
It can be shown that the schematic category has finite fibered products and that are preserved by both the forgetful to $\mathbf{CRing}\data$ and by $\Spec$. 
\begin{rem}
Heuristically, the restriction maps of $X$ being flat epimorphisms implies that the information in $X$ can be recovered from $\Spec(X)$. The other schematicity conditions can be shown to be equivalent to the existence of a certain map $\pi_X\colon \Spec(X)\to X$, i.e. to $X$ being essentially a "finite model" of $\Spec(X)$ in a topological sense.
\end{rem}
\begin{defn}
A morphism $f\colon X\to Y$ is said to be a \textit{qc-isomorphism }if $\Spec(f)$ is an isomorphism. 
\end{defn}
The class of qc-isomorphisms is a multiplicative system of arrows in $\mathbf{SchFin}$ that is maximal by definition, so the corresponding localization---Verdier quotient---defines a faithful---but not full---functor
\begin{align*}
\Spec\colon \mathbf{SchFin}_{\qc}\to \mathbf{LRS}.
\end{align*}

To study properties $\PPP$ of schematic spaces we will ask for two requisites: 
\begin{itemize}
\item A "rigorous" requisite: that $\PPP$ factors through the localization, i.e. that any representative of the qc-isomorphism class of a space or arrow determines is the whole class verifies the property or not. In other words, $\PPP$ is \textit{geometric}.
\item A "moral" requisite: $\PPP$ can be studied in terms of finite models, without applying the $\Spec$ functor. In other words, $\PPP$ is \textit{discretizable}.
\end{itemize}

In certain cases, one can "rigidify" poorly-behaved properties by studying them on certain reflective subcategories of $\mathbf{SchFin}$ that induce equivalences after localizing by qc-isomorphisms. As an example of this, see the discussion about connectedness in \cite{paper grupo etale}. In this paper we will tacitly assume that all our definitions work in this nice way, but let it be known that more technical considerations are needed for a full exposition---and that is the reason why we employ quotation marks so often to highlight seemingly ordinary notions---.
\medskip

Finally, the main result in \cite{paper grupo etale}---which, with enough work, can be written in much more geometric and elegant terms that the ones presented there---is concerned with the existence of a Galois category of "finite étale covers" for any "connected" schematic space $X$ that, when $X$ models a qc-qs scheme $S$, is naturally equivalent to the homonym Galois category of $S$. 
\begin{thm}\cite{paper grupo etale}\label{theorem galois} 
For a schematic space $X$ and a schematic morphism $\overline{\x}\colon \Spec(\Omega)\to X$ with $\Omega$ an algebraically closed field---a geometric point---, there exists a category $\Rev(X)$ and a functor $\Fib_{\overline{\x}}\colon \Rev(X)\to \mathbf{FinSet}$ such that, when $X$ is "connected", the pair $(\Rev(X), \Fib_{\overline{\x}})$ is a Galois category. We denote its fundamental group by $\pi_1^{\mathrm{et}}(X, \overline{x})$. If $S=\Spec(X)$ is a scheme, the $\Spec$ functor induces a equivalence of Galois categories $(\Rev(X), \Fib_{\overline{\x}})\simeq (\mathbf{FEt}(S), \Fib_{\Spec(\overline{\x})})$, where $\mathbf{FEt}(S)$ is the category of finite étale covers of $S$, and thus $\pi_1^{\mathrm{et}}(X, \overline{x})\simeq \pi_1^{\mathrm{et}}(S, \Spec(\overline{x}))$.
\end{thm}
\begin{rem}\label{remark pi is geometric}
In \cite{paper grupo etale} we also showed that qc-isomorphic spaces have equivalent Galois categories of finite étale covers: the construction is geometric.
\end{rem}

Of course, for a general schematic space $X$, one can consider the set of all its geometric points and define the \'etale fundamental (Stone) groupoid $\Pi_1^{\mathrm{et}}(X)$. In its general version, the Galois Theorem states that the fiber functors induce an equivalence of categories
\begin{align*}
\Rev(X)\simeq [\Pi_1^{\mathrm{et}}(X), \mathbf{FinSet}]\equiv  \Pi_1^{\mathrm{et}}(X)\text{-}\mathbf{FinSet}
\end{align*}
where the action of this groupoid is continuous. As always, this is just a particular case of more general topos-theoretic results. 

\subsection{The topology of flat immersions}\label{section topos flat immersions}
We begin by introducing the natural (pre)topology on $\mathbf{SchFin}$.

\begin{defn}
Let $f\colon X\to Y$ be a schematic morphism. We say that $f$ is \textit{flat} if $f^\sharp_x\colon\OO_{Y, f(x)}\to \OO_{X, x}$ is flat for all $x\in X$. Such flat morphism is a \textit{flat immersion} if its diagonal $\Delta_f\colon X\to X\times_YX$ is a qc-isomorphism. A flat morphism $f$ is \textit{faithfully flat} if $\Spec(f)$ is surjective.
\end{defn}

These three types of maps can be characterized in terms of the adjoint pair $(f^*, f_*)$ for quasi-coherent sheaves. We remark that a flat immersion is, by definition, a \textit{flat monomorphism} in $\mathbf{SchFin}_{\qc}$; and $\Spec(f)$ for such $f$ is a flat monomorphism of locally ringed spaces. One can show that qc-isomorphisms are exactly faithfully flat immersions.

\begin{rem}
It can be shown that a morphism $U_x\to U_y$ is a flat immersion if and only if $\OO_{Y, y}\to\OO_{X, x}$ is a flat epimorphism of rings. Since schematic spaces have flat epimorphisms of rings as restriction maps, the restriction morphisms between their basic open subsets are flat immersions---actually, between all their open subsets---. In other words, schematic spaces are colimits of (certain) affine schematic spaces via flat immersions. This class of morphisms was first shown to be important in the context of descent problems in \cite{raynaud}.
\end{rem}
\begin{lem}\label{lema inmersiones}
If $f\colon X\to Y$ is a flat immersions, $f_x^\sharp\colon\OO_{Y, f(x)}\to \OO_{X, x}$ are flat epimorphisms of rings for all $x\in X$.
\begin{proof}
They are flat by definition and the condition on the diagonal trivially translates to $\OO_{Y, f(x)}\otimes_{\OO_{X, x}}\OO_{Y, f(x)}\to \OO_{Y, f(x)}$ being an isomorphism.
\end{proof}
\end{lem}

Recall that an open immersion of schemes is a flat monomorphism (locally) of finite presentation. As such, flat immersions are like "open immersions", but without the finite presentation condition. The reader might notice the analogy with the \'etale and pro-\'etale topologies for schemes. This justifies the following notation:
\begin{defn}
Let $X$ be a schematic space. We define $X_{\mathrm{wZar}}$ to be the site of flat immersions with target $X$, whose covers are given by finite and jointly faithfully flat families of flat immersions. Similarly, we define $\mathbf{SchFin}_{\mathrm{wZar}}$ to be the "big" site of flat immersions.
\end{defn}

These sites present a number of interesting pathologies that we will describe in more detail in future papers. We shall remark a few of them:
\begin{itemize}
\item The category $X_{\mathrm{wZar}}$ is not small, only its localization $(X_\mathrm{wZar})_\qc$. Each qc-isomorphism class of open immersions is identified with a subset of $\Spec(X)$, but before localization, the collection of representatives is as large as the entire class of finite posets.
\item Since qc-isomorphisms are both flat immersions and covers, yet they are not isomorphisms, a standard descent argument shows that \textit{sheaves map qc-isomorphisms to isomorphisms}---Category theorists sometimes call morphisms with such property \textit{local isomorphisms}---. In particular, functors of points are not sheaves, because they determine spaces up to isomorphism, so the site $X_{\mathbf{wZar}}$ and its bigger analogue \textit{are not subcanonical}. However, one can show that, for any $Y\in X_{\mathrm{wZar}}$ and sheaf $\F\in X_{\mathrm{wZar}}$, there are natural bijections
\begin{align*}
\Hom_{\mathbf{PSh}(X_{\mathrm{wZar}})}&(\F, \Hom_{X_{\mathrm{wZar}}}(-, Y))\simeq\\&\simeq\Hom_{\mathbf{Sh}(X_{\mathrm{wZar}})}(\F, \Hom_{(X_{\mathrm{wZar}})_{\qc}}(-, Y))
\end{align*}
In other words, \textit{the functor of points in the localization satisfies the universal property of sheafification.}
\item We have avoided talking about sheafification in the previous points because, due to potential size issues, we cannot guarantee that such functor exists in $X_{\mathrm{wZar}}$---this is related to the inability to find bounds for refinements of covers, which may lead to pathologies, as it happens with the fpqc topology of schemes, see \cite[Theorem 5.5]{fpqc}---; but the good news is that \textit{it does exist} $(X_{\mathrm{wZar}})_\qc$. I.e. we can sheafify presheaves that factor through qc-isomorphism, which will be enough in all natural situations.
\item Endowing $X_{\mathrm{wZar}}$ with the natural sheaf of rings, it is possible to show that $\Qcoh(X_{\mathrm{wZar}})\simeq \Qcoh(X)$.
\end{itemize}

As it happens with schemes and open immersions, it is obvious that if $\X=\{X_i\}_{i\in I}$ is a diagram of schematic spaces and the transition morphisms $X_i\to X_j$ (for any $i\to j$) are flat immersions, taking $\colim_iX_i$ in the category of ringed posets yields $\Spec(\colim_iX_i)=\colim_i\Spec(X_i)$ and the resulting space is a gluing of affine schemes via flat monomorphisms of affine schemes. The problem is that, in general, it is very difficult to determine if $\colim_i X_i$ is schematic or not, due to the combinatorial nature of the definition of schematicity and the surprisingly subtle description of colimits of finite posets (see \cite[Proposition 2.4]{tesisgen posets fibra coherente}).

Our solution will be defining an object "equivalent" to $\colim_iX_i$ in the sense of representing the same locally ringed space, but whose combinatorial nature is elementary. This will be done in Sections \ref{seccion cylinder} and \ref{section schematic cylinder}. The result will be called \textit{cylinder space}, denoted $\Cyl(\X)$. 
\smallskip 

This construction is central in the theory of schematic spaces will have applications that are beyond our purposes here, but the goal for this paper is to study descent properties with respect to the topology of flat immersions. For instance, let us consider the case of the \'etale fundamental groupoid. It clearly defines a functor 
\begin{align*}
\Pi_1^{\text{\'et}}\colon\mathbf{SchFin}\to \mathbf{Gpd}_{\mathrm{Stone}}
\end{align*}
valued in the strict 2-category of Stone groupoids. Proving the Seifert-Van Kampen Theorem in its general form---for the topology of flat immersions---essentially amounts to saying that $\Pi_1^{\text{\'et}}$ maps colimits to 2-colimits. This will be the same as saying that $(\Pi_1^{\text{\'et}})^\op$ \textit{is a 2-sheaf}---thus it maps qc-isomorphisms to equivalences---. By the properties of these sites, this is equivalent to proving that it maps objects "qc-equivalent to colimits"---our cylinder spaces---to 2-colimits. However, our abstract descent result for the topology of flat immersions and cylinders will show that it is enough to prove that it is a 2-sheaf in the combinatorial topology. Such statement amounts to showing that $\Pi_1^{\text{\'et}}$ maps a very specific kind of cylinders to 2-colimits; and in some particular cases, this will even be formal.

\section{Categories of $\C\data$}

Without further ado, let $\C$ be a 1-category and $\pos$ be the category of finite posets---or arbitrary posets, being careful in that case with set-theoretic size considerations---. For a given poset $X$ and $x\in X$, let $U_x=\{x'\geq x\}$ denote the minimal open neighborhood of the point $X$. The following is well-known:
\begin{lem}\label{lemma sheaves equal diagrams}
If $\C$ has finite limits and $X\in\pos$, there is an equivalence
\begin{align*}
\mathbf{Sh}(X, \C)\simeq [X, \C]
\end{align*}
between the categories of $\C$-valued sheaves on $X$ and functors $X\to \C$. 
\begin{proof}
Each sheaf gives a functor defined by its stalks---sections at the minimal open neighborhoods---and restrictions morphisms. The converse follows from the sheaf condition and the fact that the $\{U_x\}_{x\in X}$ are a basis for the topology, so for any open $U\subseteq X$ and functor $\F\colon X\to \C$, one defines its "sections" on $U$ as $\F(U)=\lim_{x\in U}\F(x)$.
\end{proof}
\end{lem}

Now let us consider the functor to the 1-category of categories---big enough so that $\C\in\mathbf{Cat}$---
\begin{align*}
\underline{\C\data}\colon \pos&\to \mathbf{Cat}\\
X&\mapsto[X, \C]\\
f&\mapsto f^{-1}.
\end{align*}
\begin{defn}\label{definition c data}
For any $\C$, the cateory of $\C\data$ is the fibered category over $\C$ defined by the Grothendieck construction applied to the previous functor. Explicitly:
\begin{itemize}
\item $\mathrm{Ob}(\C\data)=\{\F\overset{\text{not}}{\equiv}(X, \F):X\in\pos\text{ and }\F\in[X, \C]\}$;
\item $\Hom_{\C\data}((X, \F), (Y, \G))=\{f\colon X\to Y\text{ and }f^\sharp\colon f^{-1}\G\to \F\}$; 
\item $|-|\colon \C\data \to \pos$ is the "underlying poset" structure functor.
\end{itemize}
\end{defn}
\begin{nota}
We will usually denote $\F\overset{\text{not}}{\equiv}(X, \F)$ and $X=|\F|$, unless $\C=\mathbf{CRing}$ is the category of commutative rings, in which case $\C\data$ is the category of ringed posets and we will keep the traditional notation $(X, \OO_X)$. Furthermore, for any $\F$ and $x\leq y\in |\F|$, we will denote its "restriction morphisms" by $\F_{xy}\colon  \F(x)\to \F(y)$.
\end{nota}
\begin{rem}
Note that the construction of $\C\data$ is functorial on the category: if $\Phi\colon \C\to \D$ is a functor, we have  $\Phi_*\colon \C\data\to \D\data$ induced by post-composition. 
\end{rem}

This category comes with a natural inclusion functor
\begin{align*}
i_\C\colon \C^\op\to \C\data & & c\mapsto (\star, c)
\end{align*}
analogue to the "diagonal inclusion" in categories of diagrams of a fixed shape. Due to the choice of $\star$ as the final object in $\pos$, we have the following:
\begin{lem}\label{lemma sections}
If $\C$ has finite limits (resp. colimits), the functor $i_\C$ has a left (resp. right) adjoint $\Gamma\equiv \Gamma_\C\colon \C\data\to \C^\op$ (resp. $\mathrm{L}$) called the \textit{sections} (resp. \textit{cosections}) functor. Explicitly, $\Gamma(\F)=\lim \F$ (resp. $\Gamma(\F)=\colim \F$).
\end{lem}

\begin{rem}
The terminology of Lemma \ref{lemma sections} comes from the equivalence of Lemma \ref{lemma sheaves equal diagrams}. Of course, one may assume no hypothesis on $\C$ and define sections via Yoneda at the level of $[\C^{\op}, \mathbf{Set}]\data$, only to ask if these "sheaves of sections" are representable on a case-by-case basis. One may also interpret sections via projections to the terminal poset $\pi\colon X\to \star$ by constructing $\pi_*$ right adjoint to $\pi^{-1}$.
\end{rem}

\begin{ejem}[Locally representable functors]
As a simple application of this terminology, we will give a "structured" interpretation of the concept of locally representable functor. Indeed, let $Y\colon \C\to [\C^\op, \mathbf{Set}]$ be the Yoneda embedding for $\C$ and $Y_*\colon \C\data \to [\C^\op, \mathbf{Set}]\data$ the---fully faithful---induced functor. One may think of an object in the image of $Y_*$ as a "locally representable functor". Note that, if $\C$ has finite limits, the sections of such an object are representable by the sections of the original $\C$-datum. Additionally, we shall consider the Yoneda embedding for $\C\data$, that is $Y'\colon \C\data\to [\C\data^\op, \mathbf{Set}]$. At this stage, we define a third functor
\begin{align*}
\mathcal{D}\colon [\C^\op, \mathbf{Set}]\data&\to [\C\data^\op, \mathbf{Set}]\\  \X&\mapsto\Hom_{[\C^\op, \mathbf{Set}]\data}(Y_*(-), \X)
\end{align*}
such that $\mathcal{D}\circ Y_*=Y'$---since $Y_*$ is fully faithful---. We leave as an exercise to the reader checking that $\mathcal{D}$ is fully faithful itself---recall that categories of presheaves are compactly generated by their representable functors---. In particular, if $\X$ is such that $\mathcal{D}(\X)$ is representable by some $\F\in\C\data$, one has that $Y_*(\F)\simeq\X$, in other words, "representing each $\X(p)$ by some $\F_p\in\C$ for each $p\in |\X|$ in a compatible way is equivalent to representing $\X$ by a $\C$-datum $\F$ with $\F(p)=\F_p$".
\end{ejem}

One of the main advantages of considering $\C\data$ over categories of diagrams of fixed shape is that it inherits the natural 2-categorical structure of $\pos$. More precisely, recall that $\pos$ is a strict 2-category with its 2-morphisms being, for each $X, Y\in\pos$,
\begin{align*}
\Hom_{\Hom_{\pos}(X, Y)}(f, g)=\begin{cases}\star\text{ if }f\leq g\\ \emptyset\text{ otherwise.}\end{cases}
\end{align*}

If $f, g\colon\F\to \G$ are morphisms in $\C\data$ and $|f|\leq|g|$ in $\pos$, we have a natural transformation $r_{fg}\colon f^{-1}\G\to g^{-1}\G$ given, at each $x\in |\F|$, by the restriction morphisms of $\G$. We simply ask this arrow to induce a commutative triangle, i.e. we define our 2-morphisms to be:
\begin{align*}
\Hom_{\Hom_{\C\data}(\F, \G)}(f, g)=\begin{cases}\star\text{ if }f\leq g\text{ and }g^\sharp=r_{fg}\circ f^\sharp\\ \emptyset\text{ otherwise.}\end{cases}
\end{align*}

We note that this structure generalizes the partial order defined in \cite{fernando homotopy} to study naif homotopy types of ringed posets. We also remark that $\C\data$ is actually a $\pos$-enriched category.
\medskip

It is easy to check that, if $\C$ has finite limits (resp. colimits), then $\C\data$ has finite colimits (resp. limits), described in an analogous way as in the category of ringed posets (or spaces). To approach descent problems, we are interested in computing colimits of $\C\data$, or in other words, describing the sections functor of the inclusion
\begin{align*}
i_{\C\data^{\op}}\colon \C\data\to (\C\data)^\op\data;
\end{align*}
but it turns out that we can obtain, up to a certain to-be-introduced notion of weak equivalence, a more explicit description of these colimits that does not require us to perform any 1-categorical operations on either $\pos$ or $\C$. We will call this construction the "cylinder functor". The context in which it arises naturally employs the 2-categorical structure of $\C\data$, hence, for this and other reasons, we shall devote the next section to briefly describe $\D\data$ for $\D$ a strict 2-category.

\section{The $2$-categorical case}

Let $\D$ be a strict 2-category and endow posets with the trivial 2-categorical structure. Among other possibilities, we shall consider the categories of
\begin{itemize}
\item pseudofunctors $X\to \D$ and pseudonatural transformations, $[X, \D]$;
\item pseudofunctors $X\to \D$ and lax natural transformations, $[X, \D]_{\mathrm{Lax}}$.
\end{itemize}

The Grothendieck construction for each of these possibilities now yields, as in Definition \ref{definition c data}, two different 1-categories, denoted for emphasis as $\D\data$ and $\D\data^{\mathrm{Lax}}$ respectively. In both cases, their objects are pairs $(X, \F)$ of a finite poset and a pseudofunctor, with the only difference being that a morphism $(f, f^\sharp)\colon \F\to \G$ is defined by a pseudonatural transformation $f^\sharp$ when considering it in $\D\data$ and by a Lax natural transformation when considering it in $\D\data^{\mathrm{Lax}}$.  Note that, if $\D$ is $\pos$-enriched---as is the case when $\D=\C\data$ for a 1-category $\C$---, defining such a lax natural transformation amounts to giving, for each $p\leq q\in|\F|$, 1-morphisms $\alpha_p\colon\F(p)\to \G(f(p))$ such that
\begin{align*}
\G_{f(p)f(q)}\circ \alpha_p\leq \alpha_q\circ \F_{pq},
\end{align*}
rather than asking for strict equality. 

Furthermore, in order to turn the inclusion functors
\begin{align*}
i_\D&\colon \D^\op\to \D\data,\\
i_\D^{\mathrm{Lax}}&\colon\D^\op\to \D\data^{\mathrm{Lax}},\\
&\phantom{\colon}\D\data\to \D\data^{\mathrm{Lax}}
\end{align*}
into pseudofunctors, we need to endow both categories of data with the same \textit{lax 2-categorical structure}, whose 2-morphisms are:
\begin{align*}
\Hom_{\Hom_{\D\data}(\F, \G)}(f, g)=\begin{cases}\eta\colon r_{fg}\circ g^\sharp\to f^\sharp\text{ when }|f|\leq |g|\\ \emptyset\text{ otherwise.}\end{cases}
\end{align*}
Again, if $\D$ is $\pos$-enriched, giving this lax natural transformation amounts to asking that, for $|f|\leq |g|$, we only have
\begin{align*}
r_{fg}\circ f^\sharp\leq g^\sharp.
\end{align*}
With this structure, $i_\D$ and $i_\D^{\mathrm{Lax}}$ are pseudofunctors that map any 2-morphism $\eta\colon s\to t$ in $\D$ to the 2-morphism defined by the natural transformation $\eta$, since $r_{i_\D(s)i_\D(t)}$ is the identity and the underlying posets are singletons.
\smallskip

Finally, as in the 1-categorical case, and almost by definition, we have:
\begin{prop}
The left 2-adjoint of $i_\D$ (resp. $i_\D^{\mathrm{Lax}}$) is, if it exists, the pseudolimit (resp. lax limit) of the structure pseudofunctor. We call it \textit{sections} (resp. \textit{lax sections}) functor and denote it by $\Gamma\equiv \Gamma_\D$ (resp. $\mathrm{Lax}\Gamma$).
\end{prop}

\section{The Cylinder Functor}\label{seccion cylinder}

Now we construct the lax sections functor for the 2-category $\D=\C\data^{\op}$ with $\C$ a 1-category, that is, the lax colimit functor in $\C\data$. We begin with the explicit description:
\begin{defn}
For any $\X\in(\C\data)^\op\data$, we define the \textit{cylinder of $\X$} as the $\C$-datum $\Cyl(\X)$ such that:
\begin{itemize}
\item As a set, $|\Cyl(\X)|=\coprod_{p\in|\X}|\X(p)|$. We endow it with the partial order induced by those of $|\X(p)|$ and setting that $x_p\leq y_q$---with $x_p\in|\X(p)|$ and $y_q\in|\X(q)|$---whenever $x_p\leq\X_{pq}(y_q)$.
\item The structure functor is $\Cyl(\X)(x_p)=\X(p)(x_p)$ on objects, and its restriction morphisms are given by $\X(p)_{x_px'_p}$ in each $\X(p)$ and by $$(\X_{pq})^\sharp_{y_{q}}\colon\Cyl(\X)(y_q)\to \Cyl(\X)(\X_{pq}(y_q))$$ when $p\leq q$.
\end{itemize}
\end{defn}

It is easy to check that this construction is functorial, thus we have
\begin{align*}
\Cyl\colon (\C\data)^\op\data\to \C\data.
\end{align*}
\begin{lem}\label{lema cilindros universal}
If $\C=\star$, hence $\C\data=\pos$, the functor $\Cyl$ coincides up to natural isomorphism with the lax sections functor of the inclusion $i_{\star\data^{\op}}$. In other words, $\pos$ has $\pos$-indexed lax colimits, described by $\Cyl$.
\begin{proof}
We will check that, for any $Y\in\C\data$ and $\X\in(\C\data)^{\op}\data$, there are functorial \textit{isomorphisms} of categories
\begin{align*}
\Hom_{\mathbf{pos}}(\Cyl(\X), Y)\overset{\sim}{\to}\Hom_{\pos^{\op}\data^{\mathrm{Lax}}}(\X, Y).
\end{align*}
Since $Y\equiv i_{\C\data^{\op}}(Y)$ has the terminal category $\star$ as underlying poset, there is an isomorphism $\Hom_{\pos^{\op}\data^{\mathrm{Lax}}}(\X, Y)\simeq \Hom_{[X, \pos]_{\mathrm{Lax}}}(\X, Y)$, where $X \equiv |\X|$. Now, given a morphism $f\colon\Cyl(\X)\to Y$, we have, by construction, a family of morphisms $\{f_p\colon \X(p)\to Y\}_{p\in X}$ that verify $f_p\circ \X_{pq}\leq f_q$ for $p\geq q$. This is exactly the information that defines a lax natural transformation $\X\to Y$: giving, for each $p\in X$, an arrow $\X(p)\to Y(p)=Y$ in $\pos$ and, for each $p\leq q$, a 2-morphism on the corresponding diagram, which amounts to asking that the previous inequalities hold. The converse follows from the same argument: given $g\colon\X\to Y$, the $g_p^\sharp$ are exactly the morphisms $f_p$.

Finally, saying that two morphisms $f, g\colon\Cyl(\X)\to Y$ verify $f\leq g$ is just saying that $f_p\leq g_q$ for all $p\in X$---with the previous notations---. This is precisely the notion of 2-morphism in $\pos^{\op}\data^{\mathrm{Lax}}$.
\end{proof}
\end{lem}

\begin{prop}
For any category $\C$, the functor $\Cyl$ coincides up to natural isomorphism with the lax sections functor of the inclusion $i_{\C\data^{op}}$. I.e. $\C\data$ has $\pos$-indexed lax colimits and they are described by $\Cyl$.
\begin{proof}
Again, we check that for $Y\in\C\data$ and $X\in(\C\data)^\op\data$, there are functorial \textit{isomorphisms} of categories 
\begin{align*}
\Hom_{\C\data}(\Cyl(\X), Y)\overset{\sim}{\to}\Hom_{(\C\data)^{\op}\data^{\mathrm{Lax}}}(\X, Y).
\end{align*}
The topological part of the proof has been taken care of in Lemma \ref{lema cilindros universal}, so we only need to check that such isomorphism extends to the level of $\C$-valued functors.

Given $f\colon\Cyl(\X)\to Y$, using the same notations as in the aforementioned Lemma, we have morphisms $f_p$ such that $f_\p\circ \X_{pq}\leq f_q\colon \X(q)\to Y$ topologically. This is a 2-morphism of $\C\data$ because, for each $y_q\in|\X(q)|$, 
\begin{align*}
Y_{(f_p\circ \X_{pq})(y_q)}\circ (f_p\circ \X_{pq})^\sharp_{y_q}=Y_{(f_p\circ \X_{pq})(y_q)}\circ (\X_{pq})^\sharp_{y_q}\circ (f_p)^\sharp_{\X_{pq}(y_q)}=(f_q)^\sharp_{y_q};
\end{align*}
but by the definition of $\Cyl(\X)$ and $f$, for all $p\leq q$ and $x_p=\X_{pq}(y_q)$, 
\begin{align}\label{ecuacion teorema cilindro universal}
\Cyl(\X)_{x_py_q}\circ (f_p)_{x_p}^\sharp=(f_q)_{y_q}^\sharp,
\end{align}
where $\Cyl(\X)_{x_py_q}=(\X_{pq})^\sharp_{y_q}$, as desired. The converse follows from the same relations.

At the level of morphisms, if we have arrows $f, g\colon\Cyl(\X)\to Y$ with $f\leq g$ in $\C\data$, they verify $|f|\leq |g|$ in $\pos$ and, for all $x_p\in\Cyl(\X)$,
\begin{align}\label{ecuacion teorema cilindro universal 2}
g^\sharp_{x_p}\circ Y_{f(x_p)g(x_p)}=f^\sharp_{x_p}.
\end{align}
If $\{f_q\colon \X(p)\to Y\}$ and $\{g_p\colon \X(p)\to Y\}$ are their corresponding families of morphisms in $[X, (\C\data)^{\op}]_{\mathrm{Lax}}$, there only remains to check that $f_p\leq g_p$ for all $p\in X$. Once again, $|f_p|\leq |g_p|$ by Lemma \ref{lema cilindros universal}, so we complete the proof by remarking that, for each $x_p\in|\X(p)|$, the fact that the equation \ref{ecuacion teorema cilindro universal 2} holds is equivalent to $f_p\leq g_p$ in $\C\data$.
\end{proof}
\end{prop}

Note that $\C\data$ is actually a $\pos$-enriched category, hence the universal property of $\Cyl$ is necessarily given by an \textit{isomorphism} of categories, rather than an equivalence. This means that, provided that colimits of $\C\data$ also exist, there is a natural transformation to the 1-categorical sections:
\begin{align*}
\Cyl\to \Gamma_{\C\data^{\op}}.
\end{align*}

One can make a case for this natural transformation being a "weak equivalence" relative to certain descent problems for information codified in a given collection of $\C$-datum. We will not introduce the full terminology here, since that would be a technical exercise far past our aim, but Sections \ref{section schematic} and \ref{section schematic cylinder} will put us in a particular case that hints towards this direction.

\begin{ejem}\label{example second functor}
A very important remark is that, not only $\Cyl\circ i_{\C\data^\op}$ is trivially the identity, but that every $\C$-datum is the "cilinder of its points". More precisely, for any $\C$, there is a second "obvious" inclusion functor given by post-composition with $i_\C^{\op}$:
\begin{align*}
(i_\C^{\op})_*\colon \C\data\to (\C\data)^\op\data;
\end{align*}
such that $(i_\C^{\op})_*(\F)$ has the same underlying poset as $\F$, but we "replace" each $\F(p)$ by the constant datum $(\star, \F(p))$. It is obvious that $\Cyl\circ (i_\C^{\op})_*$ is also the identity. Furthermore, there is a natural transformation
\begin{align*}
\eta_\C\colon (i_\C^{\op})_*\to i_{\C\data^\op}
\end{align*}
given by the natural projections to the terminal poset and identities in $\C$, which will be relevant when dealing with descent problems. 
\end{ejem}

\begin{prop}\label{prop fiber product}
The functor $\Cyl$ commutes with finite fibered products.
\begin{proof}
Exercise to the reader: it follows from the explicit construction. 
\end{proof}
\end{prop}
\section{The schematic cylinder}\label{section schematic cylinder}

The schematic category introduced in Section \ref{section schematic} is a non-full subcategory of $\mathbf{CRing}\data$, where $\mathbf{CRing}$ denotes the category of commutative rings with unit. In particular, the cylinder functor restricts to 
\begin{align*}
\Cyl\colon\mathbf{SchFin}^{\op}\data\to \mathbf{CRing}\data.
\end{align*}
The next few pages are devoted to characterizing $\mathbf{SchFin}^{\op}$-data whose cylinder spaces are schematic. The first justification is that such lax colimit represents up to "qc-isomorphism"---see discussion after the next Lemma---the same locally ringed space:

\begin{lem}
Given $\X\in\mathbf{SchFin}^{\op}\data$, the natural morphism of ringed spaces $\Cyl(\X)\to \Gamma(\X)$ induces an isomorphism $\Spec(\Cyl(\X))\overset{\sim}{\to}\Spec(\Gamma(\X))$. 
\begin{proof}
This follows from the fact that colimits commute with colimits. 
\end{proof}
\end{lem}

We would like to say that $\Cyl(\X)\to \Gamma(\X)$ is a qc-isomorphism, but note that we have not checked---and will not check---whether or not $\Gamma(\X)$ is schematic. However, it will be sufficient to check schematicity of $\Cyl(\X)$ for our applications---and crucial,  since we would not be able to guarantee the stability under qc-isomorphisms of the properties and constructions we are interested in dealing with otherwise---.

\begin{defn}
A ringed poset $X$ is said to be \textit{pseudo-schematic} if it has flat epimorphisms of rings as restriction maps.
\end{defn}
\begin{defn}
A ringed poset $X$ is \textit{$\mathbf{Mod}$-affine }if $\pi\colon X\to (\star, \OO_X(X))$ induces an adjoint equivalence $(\pi^*\dashv\pi_*)\colon \Qcoh(X)\to \mathbf{Mod}(\OO_X(X))$. We say that $X$ is affine if it is schematic and $\mathbf{Mod}$-affine.
\end{defn}
\begin{ejem}
Any ringed poset with a minimum $X=U_x$ is $\mathbf{Mod}$-affine.
\end{ejem}
\begin{rem}
If $X$ is pseudo-schematic, $\Qcoh(X)$ is a Grothendieck abelian category. In particular, if $X$ is also $\mathbf{Mod}$-affine, $\pi_*$ is exact.
\end{rem}
\begin{lem}\label{lema faithfully flat}
If $X$ is pseudo-schematic and $\mathbf{Mod}$-affine, the natural morphism $\OO_X(X)\to\prod_{x\in X}\OO_{X, x}$ is faithfully flat.
\begin{proof}
It suffices to see that $\coprod_{x\in X}\Spec(\OO_{X, x})\to\Spec(\OO_X(X))$ is surjective. Given a prime $\p\subseteq\OO_X(X)$ with non-zero residue field $\kappa(\p)$, the equivalence gives a non-zero module $\pi^*\kappa(\p)\neq 0$, thus there is some $x\in X$ such that $(\pi^*\kappa(\p))_x\simeq \kappa(p)\otimes_{\OO_X(X)}\OO_{X, x}\neq 0$. Geometrically, this means that the fiber of $\p$ via $\Spec(\OO_{X, x})\to \Spec(\OO_X(X))$ is non-empty, so we win.
\end{proof}
\end{lem}

\begin{defn}\label{definition qc-iso de anillados}
A morphism of ringed spaces $f\colon X\to Y$ between pseudo-schematic spaces will be called a qc-isomorphism if  $f^{-1}(U_y)$ is $\mathbf{Mod}$-affine for all $y\in Y$ and $f_\sharp\colon \OO_Y\to f_*\OO_X$ is an isomorphism.
\end{defn}
\begin{ejem}
Any ringed poset with a minimum $X=U_x$ is qc-isomorphic to $(\star, \OO_{X, x})$ via the natural projection.
\end{ejem}

In the schematic category, Definition \ref{definition qc-iso de anillados} restricts to the usual one. In this generality, we cannot even guarantee that the notion is stable under composition and base change, so the reader must think of it as an abbreviated way of storing information whose purpose will soon become clear. We would like to remark, however, that the notion of $\mathbf{Mod}$-affinity and the concept of qc-isomorphism it produces are particular cases of more abstract constructions for $\C\data$.

\begin{lem}
Given $\X\in\mathbf{SchFin}^{\op}\data$ whose restriction morphisms are flat immersions, $\Cyl(\X)$ is pseudo-schematic.
\begin{proof}
It follows from the construction, Lemma \ref{lema inmersiones} and Remark \ref{remark flat epic}.
\end{proof}
\end{lem}
Now, given $\X\in\mathbf{SchFin}^{\op}\data$ and $p\in|\X|$, denote by $U_p$ the datum induced on the open subset $U_p\subseteq|\X|$. We have \textit{qc-isomorphisms} of ringed spaces
\begin{align*}
\pi_p\colon \Cyl(U_p)\to \X(p).
\end{align*}
In general, given an open subset $U\subset|\X|$ and endowing it with the induced structure functor, we have open subsets
\begin{align*}
i_U\colon\Cyl(U)\hookrightarrow \Cyl(\X);
\end{align*}
so, for every $p, q\in |\X|$ and fixed $t\leq p, q$, we have natural morphisms
\begin{align*}
i_{pq}^p\colon\Cyl(U_p\cap U_q)\to \Cyl(U_p), & & i_{pq}^q\colon\Cyl(U_p\cap U_q)\to \Cyl(U_q);
\end{align*}
which, composing with the previous projections, induce
\begin{align*}
\pi_{pq}^t\colon \Cyl(U_p\cap U_q)\to \X(p)\times_{\X(t)}\X(q).
\end{align*}
Note that the space on the right hand side is always schematic and that, for every $(x_p, y_q)\in|\X(p)\times_{\X(t)}\X(q)|$, we have
\begin{align*}
(\pi_{pq}^t)^{-1}(U_{(x_p, y_q)})=U_{x_p}\cap U_{y_q}\subseteq |\Cyl(U_p\cap U_q)|\subseteq |\Cyl(\X)|.
\end{align*}
\begin{thm}\label{theorem schematic cylinder}
Given $\X\in\mathbf{SchFin}^\op\data$ whose restriction morphisms are flat immersions, $\Cyl(\X)$ is schematic if and only if for every $t\leq p, q$ in $|\X|$, the natural morphism $\pi_{pq}^t$ is a qc-isomorphism (\textit{a priori} of ringed posets, \textit{a posteriori} of schematic spaces).
\begin{proof}
With the technology introduced in this paper, we can only prove the "if" part, which will be the one used in our applications. Indeed, if $\pi_{pq}^t$ is a qc-isomorphism, $U_{x_p}\cap U_{y_q}$ is $\mathbf{Mod}$-affine for every $(x_p, y_q)$ as before and its global sections are isomorphic to $\OO_{\X(p), x_p}\otimes_{\OO_{\X(t), z_t}}\OO_{\X(q), y_q}$, with $z_t$ the common image of $x_p$ and $y_q$. Now, Lemma \ref{lema faithfully flat} translates exactly into the conditions of Definition \ref{definition schematic}.
\end{proof}
\end{thm}
For morphisms $f\colon\X\to \mathcal{Y}$ in $\mathbf{SchFin}^\op\data$, we can modify the previous construction to obtain, for each $p\in |\X|$ and $q\geq f(p)$, 
\begin{align*}
\rho_{pq}^f\colon \Cyl(U_p\cap f^{-1}(U_q))\to \X(p)\times_{\mathcal{Y}(f(p))}\mathcal{Y}(q).
\end{align*}

\begin{thm}\label{theorem schematic cylinder morphism}
Given a morphism $f\colon\X\to \mathcal{Y}$ in $\mathbf{SchFin}^\op\data$ and such that $\Cyl(\X)$ and $\Cyl(\mathcal{Y})$ are schematic, $\Cyl(f)$ is schematic if and only if for every $p, q\geq f(p)$, the map $\rho_{pq}^f$ is a qc-isomorphism.
\begin{proof}
We only prove the "if" part, which follows from the same results as Theorem \ref{theorem schematic cylinder} and the fact that, for $(x_p, y_q)\in\X(p)\times_{\mathcal{Y}(f(p))}\mathcal{Y}(q)$, one has $\rho_{pq}^{-1}(U_{(x_p, y_q)})=U_{x_p}\cap \Cyl(f)^{-1}(U_{y_q})$.
\end{proof}
\end{thm}
\begin{rem}
Note that, applied to a datum $\X$ with $\X(p)=(\star, A_p)$ for all $p$, Theorems \ref{theorem schematic cylinder} and \ref{theorem schematic cylinder morphism} restrict to the usual Definition of schematicity. See this in view of Example \ref{example second functor}.
\end{rem}

\begin{defn}
Given a finite family of flat immersions $\{U_i\to X\}_{i\in I}$, we define the \textit{Nerve datum} associated to it as $\U\in\mathbf{SchFin}^\op\data$ with underlying poset $|\U|=\PP^*(I)$---non-empty parts of $I$---and $\U(\Delta)=\prod_{i\in\Delta}U_i$ ---fibered product over $X$---.
\end{defn}

Note that $\U$ comes equipped with a morphism $\U\to X\equiv i_{\C\data}(X)$.
\begin{cor}\label{corollary cilinder}
If $\{U_i\to X\}$ a finite family of flat immersions, $\Cyl(\U)$ is schematic and the morphism $\Cyl(\U)\to X$ is a schematic flat immersion, which is a qc-isomorphism if and only if the family is a covering.
\begin{proof}
First, we check the condition of Theorem \ref{theorem schematic cylinder}: for $\Delta_1, \Delta_2\in|\U|$, $U_{\Delta_1}\cap U_{\Delta_2}=U_{\Delta_1\cup \Delta_2}$; but $\Cyl(U_{\Delta_1\cup\Delta_2})\to \U(\Delta_1\cup\Delta_2)$ is a qc-isomorphism, with  $\U(\Delta_1\cup\Delta_2)\simeq \U(\Delta_1)\times_{\U(\Delta_1\cap \Delta_2)}\U(\Delta_2)$ by definition. Schematicity of $\Cyl(f)$ follows from Theorem \ref{theorem schematic cylinder morphism} and a similar argument. 

The morphism $\Cyl(f)$ is flat by the local construction and its diagonal is a qc-isomorphism because, by Proposition \ref{prop fiber product}, 
\begin{align*}
\Cyl(\U)\to \Cyl(\U)\times_{X}\Cyl(\U)\simeq \Cyl(\U\times_X\U),
\end{align*}
and a morphism of $\mathbf{SchFin}^\op\data$ that is topologically the identity and a qc-isomorphism at each point, induces a qc-isomorphism between cylinder spaces (as shown by an easy computation).

Finally, $\Cyl(f)$ being faithfully flat (hence a qc-isomorphism) is clearly equivalent to $\{f_\Delta\colon U(\Delta)\to X\}_\Delta$ being a covering family, which happens if and only if the original family was a covering.
\end{proof}
\end{cor}

\section{Descent and the topos of flat immersions}\label{section descent}

Now we use the technology of the previous section to describe colimits in a sheaf-theoretic manner. In the following definition---if appropriate---, one shall consider $\mathbf{SchFin}$ as a 1-category with the \textit{trivial} 2-categorical structure.

\begin{defn}
Let $\C$ be a 1-category (resp. strict 2-category). A \textit{geometric datum} is a functor (resp. pseudofunctor) $\mathfrak{Dat}\colon \mathbf{SchFin}\to \C$ that maps qc-isomorphisms to isomorphisms (resp. equivalences); in other words, one that factors through $\mathbf{SchFin}_{\qc}$.
\end{defn}
\begin{ejem}
The functors $\Spec\colon\mathbf{SchFin}\to\mathbf{LRS}$, $\Qcoh\colon \mathbf{SchFin}\to \mathbf{Cat}^{\op}$ ---with values in the 2-category of categories---and $\Pi_1^{\text{ét}}\colon \mathbf{SchFin}\to \mathbf{Gpd}_{\mathrm{Stone}}$ ---with values in the 2-category of Stone groupoids---are all geometric data.
\end{ejem}

In the discussion that follows, let us assume that $\C$ is a 1-category; the argument also works for 2-categories, replacing isomorphisms by equivalences. In Example \ref{example second functor} we saw that there are two natural immersions of any category of $\C$-data into its category $(\C\data)^\op\data$. In this case, there is a natural transformation between functors in $[\mathbf{SchFin}, \mathbf{SchFin}_{\qc}^\op\data]$:
\begin{align*}
\eta_\C\colon (i_\C^{\op})_*\to i_{\C\data^\op}.
\end{align*}
\begin{rem}
For general ringed posets, this natural transformation is induced by the morphisms $(\star, \OO_{X, x})\to X$, which are not schematic. That is one of the reasons to consider the localized category, where it is induced by the triangles $(\star, \OO_{X, x})\leftarrow U_x\to X$.
\end{rem}
Now, given a geometric datum $\mathfrak{Dat}\colon \mathbf{SchFin}_{\qc}\to \C$, we define
\begin{align*}
\underline{\mathfrak{Dat}}:=\mathfrak{Dat}_*\circ (i_\C^\op)_* & & \mathfrak{Dat}\equiv \mathfrak{Dat}\circ i_{\C\data^\op};
\end{align*}
where $\underline{\mathfrak{Dat}}(X)$  is the $\C^\op$-datum with $|\underline{\mathfrak{Dat}}(X)|=|X|$ and structure functor $\underline{\mathfrak{Dat}}(X)(x)=\mathfrak{Dat}(\star, \OO_{X, x})$. These induce a natural transformation
\begin{align*}
\eta_{\mathfrak{Dat}}\colon \underline{\mathfrak{Dat}}\to \mathfrak{Dat}
\end{align*}
between functors in $[\mathbf{SchFin}, \C^\op\data]$, given by the projection to the point at the topological level and by the morphisms in $\C$
\begin{align*}
\mathfrak{Dat}(\star, \OO_{X, x})\overset{\sim}{\leftarrow}\mathfrak{Dat}(U_x)\to \mathfrak{Dat}(X).
\end{align*}
Composing with the sections functor $\Gamma\colon \C^\op\data\to \C$---always assuming that $\C$ has enough limits---, one arrives to the following definition:

\begin{defn}
We say that a geometric datum $\mathfrak{Dat}$ \textit{satisfies internal descent} if $\Gamma(\eta_{\mathfrak{Dat}})\colon \Gamma_*\circ \underline{\mathfrak{Dat}}\to \Gamma_*\circ \mathfrak{Dat}\equiv \mathfrak{Dat}$ is an isomorphism in $[\mathbf{SchFin}, \C]$.
\end{defn}

\begin{ejem}\label{ejemplo qcoh internal descent}
The datum $\Qcoh\colon\mathbf{SchFin}\to \mathbf{Cat}^\op$ satisfies internal descent. Indeed, since $\Qcoh(\star, \OO_{X, x})=\mathbf{Mod}(\OO_{X, x})$, this amounts to proving that the natural functor
\begin{align*}
\Qcoh(X)\to \tlim_{x\in X}\mathbf{Mod}(\OO_{X, x})
\end{align*}
is an equivalence of categories. This holds because quasi-coherent modules on ringed posets are collections of $\{M_x\}_{x\in X}$ with $M_x$ an $\OO_{X, x}$-module such that, for all $x\leq y$, the natural morphisms $M_x\otimes_{\OO_{X, x}}\OO_{X, y}\to M_y$ are isomorphisms; which coincides with the description of this pseudolimit in $\mathbf{Cat}$. The reader may notice that this result holds for arbitrary ringed posets, but that it tacitly requires the tensor-Hom adjunction for modules to hold. If one wants to extend the result to quasi-coherent sheaves of algebras, it is necessary to assume that $X$ is, at least, pseudo-schematic. This is because base changes by flat epimorphisms of rings satisfy said adjunction (left as an algebra exercise to the reader). 
\end{ejem}
\begin{prop}[External descent for nerves]\label{prop external}
If $X$ is a schematic space, $\{f_i\colon U_i\to X\}$ is a covering by flat immersions with associated nerve datum $\U$ and $\mathfrak{Dat}$ is a geometric datum satisfying internal descent, then there is a natural isomorphism
\begin{align*}
\colim_{\Delta\in|\U|}\mathfrak{Dat}(\U(\Delta))\overset{\sim}{\to}\mathfrak{Dat}(X).
\end{align*}
\begin{proof}
By Corollary \ref{corollary cilinder}, $\Cyl(\U)\to X$ is a qc-isomorphism, and since $\mathfrak{Dat}$ is geometric, one has that $\mathfrak{Dat}(\Cyl(\U))\simeq \mathfrak{Dat}(X)$. Since $\mathfrak{Dat}$ satisfies internal descent---applied at each $\U(\Delta)$---and colimits commute with colimits, 
\begin{align*}
&\mathfrak{Dat}(\Cyl(\U))\simeq \colim_{x_\Delta\in\Cyl(\U)}\mathfrak{Dat}(\star, \OO_{\U(\Delta), x_\Delta})\simeq \\&\simeq \colim_{\Delta\in|\U|}\colim_{x_\Delta\in|\U(\Delta)|}\OO_{\U(\Delta), x_\Delta})\simeq \colim_{\Delta\in|\U|}\mathfrak{Dat}(\U(\Delta)),
\end{align*}
which completes the proof.
\end{proof}
\end{prop}
\begin{ejem}
In the situation of Proposition \ref{prop external} and thanks to Example \ref{ejemplo qcoh internal descent}, we obtain that $\Qcoh(X)\simeq \tlim_{\Delta\in|\U|}\Qcoh(\U(\Delta))$. In particular, \textit{being quasi-coherent} is local in the topology of flat immersions.
\end{ejem}

\begin{thm}[External descent for topoi]\label{theorem descent}
If $\mathbf{SchFin}_\tau$ denotes the (big) site of schematic spaces with the combinatorial topology and $\mathbf{SchFin}_{\mathrm{wZar}}$ denotes the (big) site of flat immersions, the natural inclusion defines an equivalence of topoi
\begin{align*}
\mathbf{Sh}((\mathbf{SchFin}_{\qc})_\tau)\simeq \mathbf{Sh}(\mathbf{SchFin}_{\mathrm{wZar}}).
\end{align*}
Similarly, it induces equivalences between respective categories of $\C$-valued sheaves (resp. stacks) for any 1-category (resp. 2-category) that has finite poset-indexed colimits.
\begin{proof}
This is simply a reinterpretation of Proposition \ref{prop external} in terms of the language of Section \ref{section topos flat immersions}: sheaves in $\mathbf{Sh}(\mathbf{SchFin}_{\mathrm{wZar}})$ map qc-isomorphisms to isomorphisms, so they are geometric data in the sense fo this section. We shall remark that the analogous equivalence between small topoi does not hold---\textit{a priori}---because cylinders change the base space.
\end{proof}
\end{thm}
\begin{rem}
Thanks to the sheaf condition, it can be shown that a sheaf $\F$ in $\mathbf{Sh}(\mathbf{SchFin}_\tau)$ maps qc-isomorphisms to isomorphisms if and only if, for every \textit{affine} schematic space $X$, the natural morphism $\F(\star, \OO_X(X))\to F(X)$ is an isomorphism. 
\end{rem}

In other words, to prove that a presheaf in the schematic category is a sheaf in the topology of flat immersions, it is enough to see that it maps qc-isomorphisms to isomorphisms and that it is a sheaf in the combinatorial topology \textit{for every poset}. This is similar to what happens in the category of schemes for the set-theoretic topology and the Zariski site. A consequence for qc-qs schemes is the following \textbf{slogan}:
\begin{quote}
\textit{In the category of qc-qs schemes, any Zariski sheaf that can be studied through finite models \textit{is} a sheaf in the topology of flat monomorphisms of schemes and finite coverings.}
\end{quote}
Theorem \ref{theorem descent} makes the meaning of \textit{can be studied through} precise: such a sheaf $\F$ must induce a geometric datum on the schematic category that is a sheaf in the combinatorial topology. 
\section{Example: Seifert-Van Kampen Theorem}

A less trivial application comes from the étale fundamental groupoid---and group---, as promised. Let us consider the pseudofunctor
\begin{align*}
\Pi_1^{\text{ét}}\colon\mathbf{SchFin}\to \mathbf{Gpd}_{\mathrm{Stone}}
\end{align*}
to the 2-category of Stone groupoids. By Remark \ref{remark pi is geometric} it is a geometric datum, so to apply the results of the previous section it is enough to see that it is a sheaf in the combinatorial topology. This follows quite easily in two steps. Before that, we highlight that the category of finite étale covers defined without any detail in Theorem \ref{theorem galois} can be described in terms of quasi-coherent sheaves of algebras, as done in \cite{paper grupo etale}; more precisely:
\begin{align}\label{equation category}
\Rev(X)=\begin{cases}\text{"opposite category of quasi-coherent algebras }\A\\\ \text{such that }\OO_{X, x}\to \A_x\text{ is a finite étale ring map"}.\end{cases}
\end{align}
\begin{lem}\label{lemma descent of covers}
The pseudofunctor $\Rev\colon \mathbf{SchFin}_{\qc}\to \mathbf{Cat}$ defined on objects by\footnote{On 1-morphisms, we send each $f\colon X\to Y$ to the inverse image functor; since $\mathbf{SchFin}$ is considered as a 1-category, it only remains to specify invertible equivalences in $\mathbf{Cat}$ that make all suitable diagrams commute, but we can and do choose those to be the ones given by the universal property of tensor products.} $X\mapsto \Rev(X)$ satisfies internal descent.
\begin{proof}
We have to show that $\Rev(X)\simeq \tlim_{x\in X}\Rev(\star, \OO_{X, x})$. Since it is a subcategory of the category of quasi-coherent algebras, this follows from Example \ref{ejemplo qcoh internal descent}---bearing in mind the remark at the end---and the fact that the property of being finite étale at stalks is obviously local in this sense.
\end{proof}
\end{lem}
\begin{prop}\label{prop internal descen for pi1}
The pseudofunctor $\Pi_1^{\text{ét}}\colon \mathbf{SchFin}\to \mathbf{Gpd}_{\mathrm{Stone}}$ satisfies internal descent.
\begin{proof}
Following the notations of Section \ref{section descent}. By Lemma \ref{lemma descent of covers} we know that the natural transformation $\eta_{\Rev}\colon \underline{\Rev}\to \Rev$ induces an equivalence after taking sections. Composing $\Pi_1^{\text{ét}}$ with the 2-functor $\Phi\colon \mathbf{Gpd}_{\mathrm{Stone}}\to \mathbf{Cat}^{\op}$ such that $\G\mapsto \G\text{-}\mathbf{FinSet}$---with continuous action---, we obtain a commutative square of functors $[\mathbf{SchFin}, \mathbf{Cat}^{\op}]$
\begin{align*}
\xymatrix{ 
\underline{\Rev}\ar[r]\ar[d] & \Rev\ar[d]\\
\underline{\Pi_1^{\text{ét}}\text{-}\mathbf{FinSet}}\ar[r]^{\eta_{\Phi\circ \Pi_1^{\text{ét}}}} & \Pi_1^{\text{ét}}\text{-}\mathbf{FinSet};
}
\end{align*}
where the vertical arrows are isomorphisms after taking sections by the Galois Theorem for fundamental groupoids, hence $\Gamma(\eta_{\Phi\circ \Pi_1^{\text{ét}}})$ an isomorphism. 

Finally, since $\Phi$ well known to commute with pseudocolimits, one has that $\Gamma(\eta_{\Phi\circ \Pi_1^{\text{ét}}})\simeq \Phi\circ \Gamma(\eta_{\Pi_1^{\text{ét}}})$; and since this map is an equivalence and $\Phi$ is (2-)conservative by \cite[3.11]{articulo vankampen}, $\Gamma(\eta_{\Pi_1^{\text{ét}}})$ is an equivalence, which proves the statement.
\end{proof}
\end{prop}
\begin{ejem}
If a schematic space $X$ satisfies that $\Pi_1^{\text{ét}}((\star, \OO_{X, x}))=\{\star\}$---the trivial 2-category---for all $x\in X$, Proposition \ref{prop internal descen for pi1} yields
\begin{align*}
\Pi_1^{\text{ét}}(X)\simeq \tcolim_{x\in |X|}\{\star\}\simeq \widehat{\Pi_1(|X|)},
\end{align*}
where the hat denotes the profinite completion.
\end{ejem}
\begin{thm}\label{theorem van kampen}
Let $X$ be a schematic space and $\{U_i\to X\}$ be a covering by flat immersions with associated nerve datum $\U$. Then, the natural morphism
\begin{align*}
\tcolim_{\Delta\in|\U|}\Pi_1^{\text{ét}}(\U(\Delta))\to \Pi_1^{\text{ét}}(X)
\end{align*}
is an equivalence of topological groupoids. In other words, the functor $\Pi_1^{\text{ét}}$ is a (co)stack in the topology of flat immersions.
\begin{proof}
It follows from \ref{prop internal descen for pi1} and Theorem \ref{theorem descent}.
\end{proof}
\end{thm}
\begin{rem}
Note that the topological fundamental groupoid of $|\U|$ is always trivial, since any space of parts has generic point and thus is contractible to a point. One can give the statement of the Theorem in greater generality, for any $\X\in\mathbf{SchFin}^{\op}$-datum such that $\Cyl(\X)$ is schematic; and in that case the topological fundamental groupoid of $|\X|$ plays a role.
\end{rem}
\begin{cor}
If $S$ is a qc-qs scheme and $\{V_j\to S\}_{j\in J}$ is a finite cover by flat monomorphisms with associated nerve codatum $\V\colon \PP^*(J)\to \mathbf{Sch}^\op$---with $\mathbf{Sch}$ the category of schemes---, the natural morphism
\begin{align*}
\tcolim_{\Delta\in|\V|}\Pi_1^{\text{ét}}(\V(\Delta))\to \Pi_1^{\text{ét}}(S)
\end{align*}
is an equivalence of Stone groupoids, i.e. the étale fundamental groupoid of schemes is a costack in the topology of flat monomorphisms and finite covers.
\end{cor}

Finally, we can very easily specialize this result to fundamental groups, which a formulation that we deem more natural than that of \cite{stix}.
\begin{defn}
Given a schematic space $X$ and a cover by flat immersions with associated nerve datum $\U$ extended to $\PP(I)$ by $\U(\emptyset)=X$, a \textit{system of base points} $\overline{\x}_\star$ is an object
\begin{align*}
\overline{\x}_\star\in \mathrm{Ob}(\tlim_{\Delta\in|\U|}(\Pi_1^{\text{ét}}(\U(\Delta)))).
\end{align*}
\end{defn}
In other words: $\overline{\x}_\star$ is given by a collection geometric points $\overline{\x}_\Delta$ of $\U(\Delta)$ for each $\Delta$ and a collection of \textit{Tannaka paths}
$$\varphi_{\Delta\Delta'}\colon \Fib_{\overline{\x}_\Delta}\circ \Rev(X)(\Delta\to \Delta')\overset{\sim}{\to} \Fib_{\overline{\x}_{\Delta'}}$$
for each $\Delta\leq \Delta'$. Let us denote by $\overline{\x}=\overline{\x}_\emptyset$ the geometric point of $X$ given by this collection.
\begin{thm}\label{theorem ven kampen for groups}
Let $X$ be schematic and connected, $\U$ the nerve codatum associated to some covering by flat immersions such that $\U(\Delta)$ is connected, and $\overline{\x}_\star$ a system of base points. Then there is an isomorphism of topological groups
\begin{align*}
\colim_{\Delta\in|\U|}\pi_1^{\text{ét}}(\U(\Delta), \overline{\x}_\Delta)\overset{\sim}{\to}\pi_1^{\text{ét}}(X, \x)
\end{align*}
induced by conjugation the $\varphi_{\Delta\Delta'}$.
\begin{proof}
Since $X$ is connected, the natural inclusion $\pi_1^{\text{ét}}(X, \x)\to \Pi_1^{\text{ét}}(X)$ is an equivalence. Let $\mathbf{Gr}_{\text{Stone}}\subseteq \mathbf{Gpd}_{\mathrm{Stone}}$ be the category of profinite groups, which one may think set-theoretically or as $\mathbf{Top}$-enriched categories. Define the datum
\begin{align*}
\pi_1^{\text{ét}}(-, \overline{\x}_\star)\colon|\U|\to \mathbf{Gr}_{\text{Stone}}^{\op} & & \Delta\to \pi_1^{\text{ét}}(\U(\Delta), \overline{\x}_{\Delta})
\end{align*}
whose restriction morphisms given by conjugation with the $\varphi_{\Delta\Delta'}$. Since $\U(\Delta)$ is connected for every $\Delta$, the natural transformation
\begin{align*}
\pi_1^{\text{ét}}(-, \overline{\x}_\star)\to \Pi_1^{\text{ét}}\circ \U
\end{align*}
is an isomorphism of $\mathbf{Gpd}_{\mathrm{Stone}}$-valued pseudofunctors, hence it induces an isomorphism after taking sections. From this fact and Theorem \ref{theorem van kampen}, there are equivalences
\begin{align*}
\tcolim_{\Delta\in|\U|}\pi_1^{\text{ét}}(\U(\Delta), \overline{\x}_{\Delta})\overset{\sim}{\to}\tcolim_{\Delta\in|\U|}\Pi_1^{\text{ét}}(\U(\Delta))\simeq \Pi_1^{\text{ét}}(X);
\end{align*}
where the first groupoid is identified with the 1-colimit of abstract profinite groups $\colim_{\Delta\in|\U|}\pi_1^{\text{ét}}(\U(\Delta), \overline{\x}_{\Delta})$ and the last one is equivalent to $\pi_1^{\text{ét}}(X, \x)$ as remarked before. Since any equivalence between one-object categories is an isomorphism, the proof ends.
\end{proof}
\end{thm}
\begin{rem}Note that the topological Seifert-Van Kampen Theorem can be written in terms of $\C\data$: if $S$ is a quasi-compact topological space and 
\begin{align*}
\pi\colon S\to X
\end{align*}
is a finite model, we can turn $X$ into a $\mathbf{Top}^\op$-datum---with $\mathbf{Top}$ being the category of topological spaces---by setting that $X(x)=\pi^{-1}(U_x)$. If each one of these fibers is simply connected and we assume connectedness, the result recovers the classical one of McCord for $\pi_1$. For the higher homotopy groups, we are positive that should be a consequence of a Seifert-Van Kampen Theorem for fundamental homotopy groupoids thought as strict $n$-categories.
\end{rem}

\end{document}